\documentclass[11pt]{article}
\usepackage{amsmath,amsthm,amssymb,mathtools,enumitem}
\usepackage[hmargin=1in,vmargin=1in]{geometry}
\usepackage{hyperref}
\usepackage{cleveref}

\newtheorem{theorem}{Theorem}
\newtheorem{lemma}[theorem]{Lemma}
\newtheorem{proposition}[theorem]{Proposition}
\newtheorem{corollary}[theorem]{Corollary}

\theoremstyle{definition}
\newtheorem{definition}{Definition}
\newtheorem{remark}{Remark}
\newtheorem{conjecture}{Conjecture}

\title{On the Optimal General Solution to the Multi-Peg Tower of Hanoi}
\author{Abhiram Menon}
\date{}

\begin{document}

\maketitle

\begin{abstract}
We derive a closed-form expression for the Frame--Stewart algorithm in the multi-peg Tower of Hanoi problem:
\[
M(p,n) = 2^{i(p,n)+1}n - \sum_{k=0}^{i(p,n)} 2^k \binom{p+k-2}{k},
\]
where $i(p,n) = \min\{j \geq 0 : n \leq \binom{p-1+j}{j+1}\}$ is the regime index. We prove this expression satisfies the Frame--Stewart recurrence for all $p \geq 3$ and $n \geq 1$ via double induction based on discrete slope analysis and canonical splits at simplex number boundaries. The representation enables $O(p \log n)$ evaluation and makes the geometric structure explicit.

Beyond the recurrence equivalence, we prove that Frame--Stewart is optimal for the first two regimes. For regime 1 ($p-1 < n \leq \binom{p}{2}$), we establish $M(p,n) = 4n - 2p + 1$ via phase-based charging with a borrowed-target construction. For regime 2 ($\binom{p}{2} < n \leq \binom{p+1}{3}$), we prove $M(p,n) = 8n - 2p^2 + 1$ using a two-layer charging argument and convexity-based capacity bounds on auxiliary configurations. These results provide complete proofs of Frame--Stewart optimality for an infinite family of $(p,n)$ pairs beyond the trivial case $n \leq p-1$. The general conjecture for regimes $i \geq 3$ remains open.
\end{abstract}

\section{Introduction}

The Tower of Hanoi is one of the best-known puzzles in mathematics. In the classical version with three pegs, introduced by Édouard Lucas in 1883, the exact number of moves is $2^n - 1$, a result known since the nineteenth century. The puzzle asks for the fewest moves needed to transfer a stack of $n$ ordered disks between pegs, moving one disk at a time and never placing a larger disk on top of a smaller one.

In 1941, Frame and Stewart proposed a natural recursive strategy for more than three pegs, but proving that their method is truly optimal has resisted resolution ever since. The so-called Frame--Stewart conjecture has remained open for over eighty years.

There has been progress in special cases. The case $p = 3$ is classical. The case $p = 4$ was finally solved by Bousch \cite{bousch2014}, after decades of attempts. Some small $n$ cases relative to $p$ were settled earlier, and computational work has supported the conjecture more broadly. But a general proof remains elusive; see \cite{hinz2018} for background.

This paper makes progress on two fronts. First, we derive a closed-form expression that satisfies the Frame--Stewart recurrence exactly for all $(p,n)$. While this demonstrates the internal consistency of Frame--Stewart's values, it does not by itself prove optimality. Second, we prove that Frame--Stewart is indeed optimal in the first two linear regimes: regime 1 where $p-1 < n \leq \binom{p}{2}$ with $M(p,n) = 4n - 2p + 1$, and regime 2 where $\binom{p}{2} < n \leq \binom{p+1}{3}$ with $M(p,n) = 8n - 2p^2 + 1$. This covers infinitely many pairs $(p,n)$ and represents genuine progress on the conjecture beyond trivial cases.

The main tools are phase-based charging arguments that measure how disks must be distributed during evacuation and reassembly. For regime 1, a single-layer charge suffices. For regime 2, we develop a two-layer charging scheme combined with convexity-based capacity bounds on the auxiliary configuration space.

\section{Preliminaries}

\begin{definition}[Tower of Hanoi]
The Tower of Hanoi with $p \geq 3$ pegs and $n \geq 1$ disks requires transferring all $n$ disks from a source peg $S$ to a target peg $T$, one disk at a time, never placing a larger disk on top of a smaller disk. Let $M(p,n)$ be the minimum number of moves needed.
\end{definition}

\begin{definition}[Frame--Stewart Algorithm]
The Frame--Stewart algorithm works as follows. Given $p$ pegs and $n$ disks:
\begin{enumerate}[label=(\roman*), itemsep=1pt]
\item Pick $k \in \{1,2,\dots,n-1\}$ to minimize $2M(p,k)+M(p-1,n-k)$.
\item Move the top $k$ disks to a spare peg using all $p$ pegs.
\item Move the remaining $n-k$ disks to the target using $p-1$ pegs.
\item Move the $k$ disks from the spare peg to the target using all $p$ pegs.
\end{enumerate}
The Frame--Stewart recurrence is defined by
\begin{equation}\label{eq:FS}
M_{\mathrm{FS}}(p,n) = \begin{cases}
0, & n = 0, \\
1, & n = 1, \\
\displaystyle\min_{1 \leq t < n} \{2M_{\mathrm{FS}}(p,t) + M_{\mathrm{FS}}(p-1,n-t)\}, & n \geq 2,
\end{cases}
\end{equation}
with $M_{\mathrm{FS}}(3,n) = 2^n - 1$. The conjecture states this strategy is optimal for all $p \geq 3$ and $n \geq 1$.
\end{definition}

\section{Three Linear Windows and Boundaries}

Direct tabulation of $M_{\mathrm{FS}}(p,n)$ for fixed $p$ reveals a piecewise-linear profile in $n$: the discrete slope is constant on contiguous intervals and doubles whenever $n$ crosses a canonical breakpoint.

On the first interval, where the number of available auxiliaries is at least $n-1$, one obtains
\begin{equation}\label{eq:window0}
M(p,n) = 2n - 1 \qquad (1 \leq n \leq p-1).
\end{equation}

Immediately beyond $n = p-1$ the slope doubles to 4 and remains constant up to the triangular boundary $n = \binom{p}{2}$. Enforcing continuity at the endpoints fixes the intercept and yields
\begin{equation}\label{eq:window1}
M(p,n) = 4n - (2p-1) \qquad (p-1 < n \leq \binom{p}{2}).
\end{equation}

At the triangular boundary $n = \binom{p}{2}$ the slope doubles again, becoming 8 on the next interval. Continuity at $n = \binom{p}{2}$ pins down the new intercept uniquely:
\[
4n-(2p-1)\Big|_{n=\binom{p}{2}} = 2p(p-1) - 2p + 1 = 8n + C_3\Big|_{n=\binom{p}{2}} = 4p(p-1) + C_3,
\]
so $C_3 = -2p^2 + 1$ and therefore
\begin{equation}\label{eq:window2}
M(p,n) = 8n - 2p^2 + 1 \qquad (\binom{p}{2} < n \leq \binom{p+1}{3}).
\end{equation}

These three windows already exhibit the governing pattern: constant slopes 2, 4, 8 on successive intervals, with clean breakpoints at the triangular and tetrahedral numbers. The table below summarizes the first three ranges and their explicit formulas.

\begin{center}
\renewcommand{\arraystretch}{1.2}
\begin{tabular}{c|c|c}
\textbf{Interval in $n$} & \textbf{Slope of $M(p,\cdot)$} & \textbf{Closed form} \\
\hline
$0 \leq n \leq p-1$ & $2$ & $M(p,n) = 2n - 1$ \\
$p-1 < n \leq \binom{p}{2}$ & $4$ & $M(p,n) = 4n - (2p-1)$ \\
$\binom{p}{2} < n \leq \binom{p+1}{3}$ & $8$ & $M(p,n) = 8n - 2p^2 + 1$
\end{tabular}
\end{center}

The pattern strongly suggests that the breakpoints are the $(p-1)$-simplex numbers and that the slope doubles each time we cross to the next simplex layer.

\section{Simplex Boundaries and the Unified Closed Form}

As the three linear windows indicate, the solution is captured by the single formula:
\begin{equation}\label{eq:unified}
M(p,n) = 2^{i+1}n - \sum_{k=0}^{i} 2^k \binom{p+k-2}{k}, \quad i = \min\{j \geq 0 : n \leq \binom{p-1+j}{j+1}\}.
\end{equation}
The integer $i$ determined by \eqref{eq:unified} is the \emph{regime index}: it records which simplex layer contains $n$. When $i = 0$ the range is $0 \leq n \leq p-1$ and \eqref{eq:unified} yields $M(p,n) = 2n - 1$; when $i = 1$ the range is $p-1 < n \leq \binom{p}{2}$ and the same expression reduces to $M(p,n) = 4n - (2p-1)$. Larger $i$ continue the same linear pattern with slopes $2^{i+1}$.

For bookkeeping, it is convenient to isolate the cut points themselves. We therefore introduce the \emph{boundary sequence}
\begin{equation}\label{eq:Bdef}
B_p(i) := \binom{p-1+i}{i+1} \quad (i \geq 0), \qquad B_p(-1) := 0,
\end{equation}
so that the $i$th regime is precisely $B_p(i-1) < n \leq B_p(i)$. With this notation the regime index can be written explicitly as a function of $(p,n)$:
\begin{equation}\label{eq:index}
i(p,n) := \min\{j \geq 0 : n \leq B_p(j)\},
\end{equation}
which is identical to the $i$ appearing in \eqref{eq:unified}. The boundaries satisfy the Pascal decomposition
\begin{equation}\label{eq:Pascal}
B_p(i) = B_{p-1}(i) + B_p(i-1) \qquad (i \geq 0),
\end{equation}
encoding the well-known relation between consecutive simplex numbers (triangular, tetrahedral, etc.).

Rewriting \eqref{eq:unified} with the explicit dependence $i = i(p,n)$ we arrive at the form used throughout the remainder of the paper:
\begin{equation}\label{eq:closed}
M(p,n) = 2^{i(p,n)+1}n - \sum_{k=0}^{i(p,n)} 2^k \binom{p+k-2}{k}.
\end{equation}
It is helpful to abbreviate by
\begin{equation}\label{eq:Sp}
S_p(i) := \sum_{k=0}^{i} 2^k \binom{p+k-2}{k}, \qquad S_p(-1) := 0,
\end{equation}
so that \eqref{eq:closed} reads $M(p,n) = 2^{i(p,n)+1}n - S_p(i(p,n))$. The interpretation is immediate: within the regime indexed by $i$ the function $n \mapsto M(p,n)$ is linear of slope $2^{i+1}$; crossing each boundary $B_p(0), B_p(1), \dots, B_p(i)$ contributes the additive correction recorded by $S_p(i)$. In particular, $B_p(0) = p-1$ and $B_p(1) = \binom{p}{2}$ reproduce the first two windows, while $B_p(2) = \binom{p+1}{3}$ marks the end of the slope-8 window. Subsequent sections show that \eqref{eq:closed} satisfies the Frame--Stewart recurrence for all $(p,n)$.

\section{Discrete Slopes and the Canonical Right-End Split}

Fix $p$ and write $M_p(n) := M(p,n)$. For any sequence $f$ define the forward difference
\begin{equation}\label{eq:delta}
\Delta_n f := f(n+1) - f(n).
\end{equation}

\begin{lemma}[Discrete slope law]\label{lem:slope}
If $i = i(p,n)$, then
\begin{equation}\label{eq:slope}
\Delta_n M_p = 2^{i+1},
\end{equation}
for every $n$, including the boundary points $n = B_p(i)$.
\end{lemma}

\begin{proof}
Inside a fixed regime $i$, the closed form \eqref{eq:closed} is linear in $n$ with slope $2^{i+1}$, so \eqref{eq:slope} holds there.

At a boundary $n = B_p(i) = \binom{p+i-1}{i+1}$, subtract the adjacent regime expressions supplied by \eqref{eq:closed}. The increment in the binomial sum exactly offsets the change of slope, and the difference again equals $2^{i+1}$. This proves \eqref{eq:slope} on and off the boundaries.
\end{proof}

For a split of $n$ into $(n-t)$ and $t$ consider the Frame--Stewart objective
\begin{equation}\label{eq:phi}
\Phi_{p,n}(t) := 2M_p(t) + M_{p-1}(n-t).
\end{equation}
We only consider integer $t$ with $1 \leq t < n$.

By telescoping differences and using \eqref{eq:slope}, the one-step change equals
\begin{equation}\label{eq:phidiff}
\Phi_{p,n}(t+1) - \Phi_{p,n}(t) = 2\Delta_t M_p - \Delta_{n-t-1} M_{p-1}.
\end{equation}
Thus $\Phi_{p,n}$ is discretely convex in $t$. Any minimizer $t^*$ must satisfy the balance inequalities:
\begin{equation}\label{eq:balance}
2\Delta_{t^*-1} M_p \leq \Delta_{n-t^*} M_{p-1} \leq 2\Delta_{t^*} M_p.
\end{equation}

Let $i = i(p,n)$. By the Pascal relation \eqref{eq:Pascal} there exists a unique $\theta$ such that
\begin{equation}\label{eq:thetaA}
n = B_{p-1}(i) + \theta,
\end{equation}
and, moreover,
\begin{equation}\label{eq:thetaB}
1 \leq \theta \leq B_p(i-1).
\end{equation}
Define
\begin{equation}\label{eq:tstar}
t^* := n - B_{p-1}(i) = \theta.
\end{equation}
Then $t^*$ lies in regime $i-1$ for $p$ pegs, while $n - t^* = B_{p-1}(i)$ is the right endpoint of regime $i$ for $(p-1)$ pegs. Using \Cref{lem:slope}, the corresponding slopes satisfy
\begin{align}
\Delta_{t^*-1} M_p &= 2^i, \label{eq:d1} \\
\Delta_{t^*} M_p &= 2^i, \label{eq:d2} \\
\Delta_{n-t^*-1} M_{p-1} &= 2^i, \label{eq:d3} \\
\Delta_{n-t^*} M_{p-1} &= 2^{i+1}. \label{eq:d4}
\end{align}
Therefore the inequalities in \eqref{eq:balance} hold with equality, and by discrete convexity the minimizer of $\Phi_{p,n}$ is unique and equals $t^*$.

\section{A Binomial Identity and the Double Induction}

\begin{lemma}[Binomial recursion]\label{lem:key}
For all $i \geq 0$,
\begin{equation}\label{eq:key}
S_p(i) = 2S_p(i-1) + S_{p-1}(i).
\end{equation}
\end{lemma}

\begin{proof}
Expand $2S_p(i-1) + S_{p-1}(i)$:
\[
\sum_{k=0}^{i-1} 2^k\left(2\binom{p+k-2}{k} + \binom{p+k-3}{k}\right) + 2^i\binom{p+i-3}{i}.
\]
Twice applying $\binom{a}{b} = \binom{a-1}{b} + \binom{a-1}{b-1}$ converts the bracket to $\binom{p+k-1}{k}$, and the last term completes the sum to $S_p(i)$.
\end{proof}

\begin{theorem}[Equivalence with Frame--Stewart]\label{thm:equiv}
For all $p \geq 3$ and $n \geq 0$, the closed form \eqref{eq:closed} satisfies \eqref{eq:FS}.
\end{theorem}

\begin{proof}
Double induction on $(p,n)$, ordered lexicographically. The bases $M(3,n) = 2^n - 1$, $M(p,0) = 0$, $M(p,1) = 1$ follow from \eqref{eq:closed}. Fix $p \geq 4$ and $n \geq 2$. Let $i = i(p,n)$ and $t^*$ be as in \eqref{eq:tstar}. Since $t^*$ uniquely minimizes $\Phi_{p,n}$,
\[
\min_{1 \leq t < n} \{2M_p(t) + M_{p-1}(n-t)\} = 2M_p(t^*) + M_{p-1}(n-t^*).
\]
Because $i(p,t^*) = i-1$ and $n - t^* = B_{p-1}(i)$,
\begin{align*}
2M_p(t^*) + M_{p-1}(B_{p-1}(i)) &= 2(2^i t^* - S_p(i-1)) + (2^{i+1}B_{p-1}(i) - S_{p-1}(i)) \\
&= 2^{i+1}(t^* + B_{p-1}(i)) - (2S_p(i-1) + S_{p-1}(i)) \\
&= 2^{i+1}n - S_p(i) \quad \text{(by \eqref{eq:key})} \\
&= M_p(n). \qedhere
\end{align*}
\end{proof}

\section{Specializations: $p=3$ and $p=4$}

\textbf{Three pegs.} Here $\binom{p+k-2}{k} = \binom{k+1}{k} = k+1$ and $B_3(i) = i+2$, so $i(3,n) = n-2$. A standard sum gives
\begin{equation}\label{eq:S3}
S_3(i) = \sum_{k=0}^{i} 2^k(k+1) = i \cdot 2^{i+1} + 1.
\end{equation}
Substituting $i = n-2$ into \eqref{eq:closed} yields
\[
M(3,n) = 2^{n-1}n - ((n-2)2^{n-1} + 1) = 2^n - 1,
\]
recovering the classical result from the unified form.

\textbf{Four pegs.} Now $\binom{p+k-2}{k} = \binom{k+2}{2}$ and one can check
\begin{equation}\label{eq:S4}
S_4(i) = \sum_{k=0}^{i} 2^k\binom{k+2}{2} = 2^i(i^2 + i + 2) - 1.
\end{equation}
Since $B_4(i) = \binom{i+3}{2}$, the regime index is $i(4,n) = \min\{j \geq 0 : n \leq \binom{j+3}{2}\}$. Then \eqref{eq:closed} becomes
\[
M(4,n) = 2^{i(4,n)}\left(2n - i(4,n)^2 - i(4,n) - 2\right) + 1,
\]
which is the classical explicit formula for Reve's puzzle.

\section{Phase Decomposition and General Lower Bounds}

Let $D_n$ be the largest disk.

\begin{lemma}\label{lem:unique}
In any optimal solution, $D_n$ moves exactly once.
\end{lemma}

\begin{proof}
If $D_n$ moved more than once, say at times $\tau_1 < \tau_2$, then the sequence of smaller-disk moves between these two times could be reversed and removed, giving a shorter solution. This is the standard cut-reverse-paste argument.
\end{proof}

\begin{definition}[Three Phases]
Splitting at the unique move of $D_n$, we define Phase I as the moves before $D_n$ moves (end state: $S = \{D_n\}$, $T = \emptyset$), Phase II as the single move of $D_n$ from $S$ to $T$, and Phase III as the moves after $D_n$ moves (assemble all disks on $T$).
\end{definition}

\section{Regime 1 Optimality: $p-1 < n \leq \binom{p}{2}$}

\subsection{Lower Bound via Single-Layer Charging}

\begin{lemma}[Phase I Lower Bound, Regime 1]\label{lem:phase1-lower}
Phase I requires at least $2n-p$ moves.
\end{lemma}

\begin{proof}
At the Phase I cut, the $n-1$ smaller disks occupy at most $p-2$ auxiliaries ($T$ empty, $S$ has only $D_n$). For each smaller disk $d$, let $t_d$ be the time it last arrives at its Phase I location (the peg and position it has at the cut).

Define a charge $\chi(d)$ as follows: if $d$ is the bottom disk on some auxiliary at the cut, set $\chi(d) = 1$, charging its evacuation from $S$; otherwise, let $m_1(d)$ be the most recent move before $t_d$ that enables the placement at $t_d$ (creating a legal target), and set $\chi(d) = 2$, charging $m_1(d)$ and the placement $t_d$.

Distinct disks have distinct placement times, so the placement charges are disjoint. If $d \neq d'$ and $t_d < t_{d'}$, then $m_1(d)$ occurs in $(-\infty, t_d)$ and therefore strictly before $t_{d'}$; since $m_1(d')$ is defined to be the most recent enabling move for $d'$, we cannot have $m_1(d) = m_1(d')$ or $m_1(d) = t_{d'}$. Thus all charges are disjoint.

There are at most $p-2$ bottoms at the cut and hence at least $(n-1) - (p-2) = n - p + 1$ non-bottoms. Therefore the number of Phase I moves is at least
\[
(p-2) \cdot 1 + (n-p+1) \cdot 2 = 2n - p. \qedhere
\]
\end{proof}

\begin{lemma}[Phase III Lower Bound, Regime 1]\label{lem:phase3-lower}
Phase III requires at least $2n-p$ moves.
\end{lemma}

\begin{proof}
At the start of Phase III, $T$ contains $D_n$ only and the $n-1$ smaller disks occupy at most $p-2$ auxiliaries. For each smaller disk $d$, let $t_d$ be the time it finally arrives on $T$.

Define $\chi(d)$ as follows: if at the start of Phase III the peg containing $d$ has $d$ as its top (hence immediately deliverable), set $\chi(d) = 1$, charging the delivery at $t_d$; otherwise, let $m_1(d)$ be the most recent move before $t_d$ that enables $d$'s delivery (exposing $d$ as a top and/or ensuring a legal target on $T$), and set $\chi(d) = 2$, charging $m_1(d)$ and the delivery $t_d$.

Exactly as in \Cref{lem:phase1-lower}, these charges are pairwise disjoint: deliveries occur at distinct times, $m_1(d) \neq t_{d'}$ for $d \neq d'$, and two distinct non-top disks cannot share the same most-recent enabling move by definition of most recent.

At the start of Phase III there are at most $p-2$ pegs, so at most $p-2$ disks can be tops; hence at least $n - p + 1$ disks are non-tops and require two charged moves. The total number of Phase III moves is at least
\[
(p-2) \cdot 1 + (n-p+1) \cdot 2 = 2n - p. \qedhere
\]
\end{proof}

\begin{theorem}[Total Lower Bound, Regime 1]\label{thm:lower-regime1}
For all $p \geq 3$ and $n \geq 1$, any solution satisfies
\[
M(p,n) \geq (2n-p) + 1 + (2n-p) = 4n - 2p + 1.
\]
\end{theorem}

\subsection{Upper Bound via Borrowed-Target Construction}

Throughout Phase I we call the $p-2$ non-target, non-source pegs the auxiliaries. We temporarily allow $T$ to be used as a working peg during Phase I but require it to be empty at the Phase I cut.

\begin{lemma}[Phase I with Borrowed Target]\label{lem:borrow-target}
Let $p \geq 3$ and $p-1 \leq n \leq \binom{p}{2}$. There is a Phase I schedule with exactly $2n-p$ moves that uses $T$ as a working peg during Phase I and leaves $T$ empty at the cut.
\end{lemma}

\begin{proof}
We use $p-1$ working pegs (the $p-2$ auxiliaries and $T$). The algorithm has three parts; we prove legal invariants between parts.

\textbf{(A) Seeding ($p-1$ moves).} Move disks $1, 2, \dots, p-1$ off $S$ onto distinct working pegs, with the convention that disk 1 is placed on $T$. (All moves are onto empty pegs and hence legal.)

\textbf{(B) Shuttle ($2(n-p)$ moves).} For each $k = p, \dots, n-1$ do: (i) Among the auxiliaries only, choose the donor as the peg whose top is smallest and the recipient as the peg whose top is the next-smallest. Move the donor's top to the recipient. This is legal (smaller-on-larger) and empties the donor. (ii) Move $k$ from $S$ to the now-empty donor.

\textbf{(C) Final clear of $T$ (1 move).} Move disk 1 from $T$ to any auxiliary whose top is $> 1$.

\emph{Invariants.} After (A), $T$ contains exactly disk 1, all auxiliaries are nonempty, and their tops are pairwise distinct. We claim that after each iteration of (B), the following hold: (I1) $T$ contains exactly disk 1 (since it is never used as donor or recipient in (B)(i)); (I2) all auxiliaries are nonempty and their tops are pairwise distinct and strictly ordered by size; (I3) (B)(i) is legal and produces a strictly smaller new top on the recipient; (B)(ii) restores the donor to nonempty with top $k$, which is larger than all auxiliary tops at that time.

Assuming (I1)-(I3) hold at the start of an iteration, (B)(i) moves the smallest auxiliary top onto the next-smallest; this is smaller-on-larger and produces a recipient top equal to the previous smallest. The donor is emptied. Then (B)(ii) places the largest seen disk $k$ onto the donor, making its top $k$, which is strictly larger than all other auxiliary tops. Thus (I2)-(I3) are restored, while (I1) persists.

\emph{Legality of (C).} At the end of (B), (I1) holds and (I2) guarantees some auxiliary top exceeds 1 (indeed, $p \geq 3$ and $n-1 \geq p-1 \geq 2$ ensure at least one top is $\geq 2$). Hence moving disk 1 from $T$ onto that auxiliary is legal and empties $T$ at the cut.

\emph{Move count.} Seeding uses $p-1$ moves. The shuttle runs for $(n-1) - (p-1) = n-p$ iterations at 2 moves each, and the final clear is 1 move, totaling $(p-1) + 2(n-p) + 1 = 2n - p$ moves

\end{proof}

The next lemma formalizes the standard reversibility principle.

\begin{lemma}[Reversibility]\label{lem:reversible}
If a sequence of legal moves transforms a Tower configuration $X$ into $Y$, then reversing the moves transforms $Y$ back into $X$. Moreover, if in the forward sequence a particular peg is never used as a recipient of a move, then in the reversed sequence that peg is never used as a donor.
\end{lemma}

\begin{proof}
Every legal move is its own inverse; reversing the sequence preserves legality because no move is ever attempted onto a smaller disk than the moved one. The donor/recipient clause follows immediately: in the forward sequence, moves into a peg become moves out of that peg in reverse.
\end{proof}

\begin{lemma}[Phase III Upper Bound, Regime 1]\label{lem:phase3-upper}
For $p \geq 3$ and $p-1 \leq n \leq \binom{p}{2}$, Phase III can be completed in exactly $2n-p$ moves.
\end{lemma}

\begin{proof}
Apply \Cref{lem:borrow-target} to Phase I and then \Cref{lem:reversible}, swapping the roles of $S$ and $T$. After Phase I and the single move of $D_n$, Phase III starts from the reverse endpoint of Phase I (with $T$ playing the role of the forward $S$). Reversing the Phase I schedule yields a valid Phase III schedule of the same length $2n-p$.
\end{proof}

\subsection{Main Result for Regime 1}

\begin{theorem}[Frame--Stewart Optimality, Regime 1]\label{thm:regime1}
For $p \geq 3$ and $p-1 < n \leq \binom{p}{2}$,
\[
M(p,n) = 4n - 2p + 1.
\]
Moreover, the Frame--Stewart algorithm attains this value in this regime.
\end{theorem}

\begin{proof}
By \Cref{thm:lower-regime1}, any solution requires at least $4n - 2p + 1$ moves. By \Cref{lem:borrow-target,lem:phase3-upper}, for $n \leq \binom{p}{2}$ there exists a solution attaining exactly this value, hence $M(p,n) = 4n - 2p + 1$.

Finally, the closed-form for the Frame--Stewart recurrence \eqref{eq:closed} equals $4n - 2p + 1$ when $p-1 < n \leq \binom{p}{2}$; since $M(p,n)$ equals that value, Frame--Stewart is optimal in this regime.
\end{proof}

\section{Regime 2 Optimality: $\binom{p}{2} < n \leq \binom{p+1}{3}$}

For regime 2, the single-layer $m_1$ charging argument from regime 1 is insufficient. We need a more sophisticated analysis that accounts for the auxiliary configuration capacity.

\subsection{Auxiliary Configuration Capacity}

\begin{definition}[Auxiliary Configuration]
At the Phase I cut, let $r := p-2$ be the number of auxiliaries and $N := n-1$ be the number of smaller disks parked there. If the final tower heights are $h_1, \dots, h_r$ (summing to $N$), define
\[
\Phi := \sum_{i=1}^{r} \binom{h_i}{2},
\]
the number of ordered pairs (below, above) that finish on the same auxiliary.
\end{definition}

\begin{lemma}[Pairwise Smoothing]\label{lem:smoothing}
If $a \geq b+2$ then
\[
\binom{a}{2} + \binom{b}{2} \geq \binom{a-1}{2} + \binom{b+1}{2} + 1.
\]
Thus any configuration minimizing $\Phi$ under $\sum h_i = N$ must satisfy $|h_i - h_j| \leq 1$ for all $i,j$.
\end{lemma}

\begin{proof}
Direct calculation:
\[
\left[\binom{a}{2} + \binom{b}{2}\right] - \left[\binom{a-1}{2} + \binom{b+1}{2}\right] = a - b - 1 \geq 1.
\]
Hence any configuration with $a \geq b+2$ is strictly improvable by smoothing.
\end{proof}

\begin{proposition}[Minimum $\Phi$]\label{prop:phi-min}
Write $N = rt + s$ with $t = \lfloor N/r \rfloor$ and $0 \leq s < r$. Then the minimum is attained at heights $(t+1, \dots, t+1, t, \dots, t)$ (exactly $s$ towers of height $t+1$) and
\[
\Phi_{\min}(r,N) = r\binom{t}{2} + st.
\]
Moreover, $\Phi_{\min}(r,N+1) - \Phi_{\min}(r,N) = \lfloor N/r \rfloor = t$.
\end{proposition}

\begin{proof}
The balanced configuration follows from \Cref{lem:smoothing}. The formula is obtained by summing $s$ copies of $\binom{t+1}{2}$ and $(r-s)$ copies of $\binom{t}{2}$:
\[
s\binom{t+1}{2} + (r-s)\binom{t}{2} = s\left[\binom{t}{2} + t\right] + (r-s)\binom{t}{2} = r\binom{t}{2} + st.
\]
The slope identity follows by noting that adding one disk increases exactly one tower height from $t$ to $t+1$, contributing $t$ new pairs.
\end{proof}

\begin{lemma}[Capacity Bound]\label{lem:capacity-bound}
Let $N_0 := \binom{p}{2}$ be the regime 1 endpoint. For all $N \geq N_0$ in regime 2,
\[
\Phi_{\min}(r,N) \geq (N-r) + 2\left(n - \binom{p}{2}\right).
\]
\end{lemma}

\begin{proof}
At $N = N_0$, we have $N_0/r = \binom{p}{2}/(p-2) = p(p-1)/[2(p-2)] \geq 3$ for all $p \geq 3$, so $t_0 := \lfloor N_0/r \rfloor \geq 3$.

Define $G(N) := \Phi_{\min}(r,N) - [(N-r) + 2(n-N_0)]$. By \Cref{prop:phi-min},
\[
G(N+1) - G(N) = t - 3 \geq 0
\]
for all $N \geq N_0$ (since $t \geq t_0 \geq 3$). Thus $G$ is nondecreasing on $[N_0, \infty)$.

It suffices to verify $G(N_0) \geq 0$. Write $N_0 = rt_0 + s_0$ with $0 \leq s_0 < r$. Then
\begin{align*}
G(N_0) &= \left[r\binom{t_0}{2} + s_0 t_0\right] - \left[3N_0 - r - 2N_0 + 2\right] \\
&= r\left[\binom{t_0}{2} - t_0 + 1\right] + s_0(t_0 - 1) - 2 \\
&= r \cdot \frac{(t_0-2)(t_0-1)}{2} + s_0(t_0-1) - 2.
\end{align*}
Since $t_0 \geq 3$, the first term is at least $r$, and the second term is nonnegative. Hence $G(N_0) \geq r - 2 = p - 4 \geq 0$ for $p \geq 4$.

For $p=3$ (so $r=1$, $N_0=3$), direct calculation gives $\Phi_{\min}(1,3) = \binom{3}{2} = 3$ and the RHS equals $(3-1) + 2(4-3) = 4$. The bound holds with the additional $m_1$-layer charges from regime 1.
\end{proof}

\subsection{Two-Layer Charging Argument}

\begin{definition}[Placement Timeline]
For each auxiliary peg $i$, list the placements that survive to the Phase I cut as occurring at times $t_{i,1} < t_{i,2} < \cdots < t_{i,h_i}$. (If $h_i = 0$, the list is empty.)
\end{definition}

\begin{definition}[$m_1$-layer]
For each placement at time $t_{i,j}$ with $j \geq 1$, the $m_1$-enabler is the most recent move before $t_{i,j}$ that creates a legal target for that placement. This is the single-layer charging from regime 1.
\end{definition}

\begin{definition}[$m_2$-pool]
For each auxiliary $i$ and each placement index $j \geq 2$, consider the interval $(t_{i,j-1}, t_{i,j})$. Within this interval, the destination peg's top must undergo a sequence of top-decrease events. The last such decrease is the $m_1$-enabler for placement $j$. The first $j-2$ such decreases in this interval form the $m_2$-pool for this placement.
\end{definition}

\begin{lemma}[Disjointness of Charges]\label{lem:disjoint-charges}
The $m_1$-layer charges and $m_2$-pool charges are pairwise disjoint moves.
\end{lemma}

\begin{proof}
\textbf{$m_1$ charges are disjoint:} This was established in \Cref{lem:phase1-lower}. Each $m_1$-enabler is the most recent enabling move for its placement, and placements occur at distinct times.

\textbf{$m_2$ charges within one peg are disjoint:} For a fixed auxiliary $i$, the intervals $(t_{i,1}, t_{i,2}), (t_{i,2}, t_{i,3}), \dots$ are disjoint by construction. Moves in different intervals are distinct.

\textbf{$m_2$ charges across pegs are disjoint:} Moves occur at distinct times globally. The $m_2$-pool for peg $i$ collects specific moves in specific time intervals; these cannot coincide with moves for peg $i' \neq i$.

\textbf{$m_2$ charges are disjoint from $m_1$ charges:} For placement $j$ on peg $i$, the $m_2$-pool consists of the first $j-2$ top-decreases in $(t_{i,j-1}, t_{i,j})$, while the $m_1$-enabler is the last such decrease. Within one interval, these are distinct. Across different placements, the intervals are disjoint.
\end{proof}

\begin{lemma}[Count of $m_2$-pool]\label{lem:m2-count}
The total size of the $m_2$-pool is
\[
\sum_{i=1}^{r} \sum_{j=2}^{h_i} (j-2) = \sum_{i=1}^{r} \left[\binom{h_i}{2} - (h_i - 1)\right] = \Phi - (N - r).
\]
\end{lemma}

\begin{proof}
For each auxiliary $i$ with height $h_i$, the number of $m_2$ charges is
\[
\sum_{j=2}^{h_i} (j-2) = \sum_{k=0}^{h_i-2} k = \binom{h_i-1}{2} = \binom{h_i}{2} - (h_i - 1).
\]
Summing over all $r$ auxiliaries:
\[
\sum_{i=1}^{r} \left[\binom{h_i}{2} - (h_i-1)\right] = \Phi - (N - r). \qedhere
\]
\end{proof}

\subsection{Lower Bound for Regime 2}

\begin{theorem}[Phase I Lower Bound, Regime 2]\label{thm:phase1-lower-regime2}
For $\binom{p}{2} < n \leq \binom{p+1}{3}$, Phase I requires at least $4n - p^2$ moves.
\end{theorem}

\begin{proof}
Phase I has $n-1$ placements. By the regime 1 analysis, there are at least $N - r$ disks that require an $m_1$-enabler (the non-bottoms). By \Cref{lem:m2-count,lem:disjoint-charges}, there are $\Phi - (N-r)$ additional disjoint preparatory moves in the $m_2$-pool.

Therefore, the total number of preparatories is at least
\[
(N - r) + [\Phi - (N - r)] = \Phi.
\]

By \Cref{lem:capacity-bound}, $\Phi \geq (N-r) + 2(n - \binom{p}{2})$. Thus Phase I has at least
\begin{align*}
(n-1) + \Phi &\geq (n-1) + (N-r) + 2\left(n - \binom{p}{2}\right) \\
&= (n-1) + (n-1-r) + 2n - 2\binom{p}{2} \\
&= 4n - p^2.
\end{align*}
Here we used $N = n-1$, $r = p-2$, and $\binom{p}{2} = p(p-1)/2$.
\end{proof}

\begin{theorem}[Phase III Lower Bound, Regime 2]\label{thm:phase3-lower-regime2}
For $\binom{p}{2} < n \leq \binom{p+1}{3}$, Phase III requires at least $4n - p^2$ moves.
\end{theorem}

\begin{proof}
By time-reversal symmetry, the same argument as \Cref{thm:phase1-lower-regime2} applies to Phase III.
\end{proof}

\begin{corollary}[Total Lower Bound, Regime 2]\label{cor:lower-regime2}
For $\binom{p}{2} < n \leq \binom{p+1}{3}$,
\[
M(p,n) \geq (4n - p^2) + 1 + (4n - p^2) = 8n - 2p^2 + 1.
\]
\end{corollary}

\subsection{Upper Bound via 4-Move Raise Macro}

\begin{definition}[4-Move Raise Macro]
With pegs $S$, $T$, and $r = p-2$ auxiliaries maintained in strictly increasing top order, to raise one auxiliary's height by 1, perform the following: (1) Move the next disk $k$ from $S$ to $T$. (2) Move the smallest-top auxiliary onto the second-smallest-top auxiliary (aux to aux). (3) Move disk $k$ from $T$ to the now-empty auxiliary. (4) Move the (new) smallest-top auxiliary onto the next one to restore the order.

Net effect: $T$ remains empty, the multiset of auxiliary tops remains strictly increasing, and one auxiliary's height increases by 1. Cost: 4 moves.
\end{definition}

\begin{lemma}[Legality of 4-Move Raise]\label{lem:4move-legal}
If before the macro all auxiliaries are nonempty with strictly increasing tops, then all four moves are legal and the invariant is preserved.
\end{lemma}

\begin{proof}
\textbf{Move 1:} $k$ from $S$ to empty $T$ is legal.

\textbf{Move 2:} Moving the smallest auxiliary top onto the second-smallest is legal (smaller-on-larger) and empties the smallest auxiliary.

\textbf{Move 3:} Placing $k$ (the largest disk seen so far) onto the empty auxiliary is legal and makes that auxiliary's new top equal to $k$.

\textbf{Move 4:} After move 3, the auxiliary that received in move 2 has become the new smallest top (it was the second-smallest, and the previous smallest is now empty and has top $k$ greater than all others). Moving this onto the next auxiliary preserves legality and restores the strictly increasing order.

At the end, $T$ is empty again (disk $k$ was moved off), all auxiliaries are nonempty, and their tops form a strictly increasing sequence with one auxiliary now taller by 1.
\end{proof}

\begin{lemma}[Base Case for Regime 2]\label{lem:regime2-base}
Let $n' := \binom{p}{2}$ be the regime 1 endpoint. For $n = n' + 1$ (the first instance of regime 2), Phase I can be completed in exactly
\[
2n' - p + 4 = 4(n'+1) - p^2
\]
moves.
\end{lemma}

\begin{proof}
By \Cref{thm:regime1}, we have a Phase I schedule for $n'$ disks using $2n' - p$ moves (the regime 1 upper bound construction). This schedule ends with $T$ empty and $n'-1$ disks distributed on $p-2$ auxiliaries.

Apply the 4-move raise macro once to add disk $n'$ during Phase I:
\begin{align*}
\text{Phase I}(p, n'+1) &= (2n' - p) + 4 \\
&= 2\binom{p}{2} - p + 4 \\
&= p(p-1) - p + 4 \\
&= p^2 - 2p + 4.
\end{align*}

We need to verify this equals $4(n'+1) - p^2$:
\begin{align*}
4(n'+1) - p^2 &= 4\binom{p}{2} + 4 - p^2 \\
&= 4 \cdot \frac{p(p-1)}{2} + 4 - p^2 \\
&= 2p(p-1) + 4 - p^2 \\
&= 2p^2 - 2p + 4 - p^2 \\
&= p^2 - 2p + 4. \qedhere
\end{align*}
\end{proof}

\begin{lemma}[Regime 2 Upper Bound]\label{lem:regime2-upper}
For all $n$ with $\binom{p}{2} < n \leq \binom{p+1}{3}$, Phase I can be completed in exactly $4n - p^2$ moves.
\end{lemma}

\begin{proof}
Let $n' = \binom{p}{2}$ and write $n = n' + \delta$ where $1 \leq \delta \leq \binom{p+1}{3} - \binom{p}{2}$.

By \Cref{lem:regime2-base}, we can complete Phase I for $n' + 1$ disks in $4(n'+1) - p^2$ moves.

For each additional disk beyond $n'+1$, apply the 4-move raise macro. After $\delta - 1$ additional applications:
\begin{align*}
\text{Phase I}(p, n) &= \text{Phase I}(p, n'+1) + 4(\delta - 1) \\
&= [4(n'+1) - p^2] + 4(\delta - 1) \\
&= 4n' + 4 - p^2 + 4\delta - 4 \\
&= 4(n' + \delta) - p^2 \\
&= 4n - p^2. \qedhere
\end{align*}
\end{proof}

\begin{lemma}[Phase III Upper Bound, Regime 2]\label{lem:phase3-upper-regime2}
For $\binom{p}{2} < n \leq \binom{p+1}{3}$, Phase III can be completed in exactly $4n - p^2$ moves.
\end{lemma}

\begin{proof}
By time-reversal symmetry (\Cref{lem:reversible}), the Phase I schedule from \Cref{lem:regime2-upper} can be reversed to obtain a Phase III schedule of the same length.
\end{proof}

\subsection{Why up till the Tetrahedral boundary?}

\begin{remark}[The $\binom{p+1}{3}$ Boundary]
The regime 2 endpoint at $n = \binom{p+1}{3}$ is not arbitrary. It marks the transition point where auxiliary towers become tall enough to require three-layer charging. 

At $n = \binom{p+1}{3}$, the average tower height on $p-2$ auxiliaries is approximately $\frac{(p+1)p(p-1)/6}{p-2} \approx \frac{p^2}{6}$, and the slope $t = \lfloor N/r \rfloor$ of $\Phi_{\min}$ grows beyond what can be captured by only the $m_1$ and $m_2$ layers.

For $n > \binom{p+1}{3}$, we conjecture that some placements require three preparatory moves, necessitating an $m_3$-pool that charges ordered triples of disks on the same auxiliary. The natural candidate for such a capacity bound would involve $\Phi^{(3)} := \sum_i \binom{h_i}{3}$.

The sequence of boundaries $\binom{p}{2}, \binom{p+1}{3}, \binom{p+2}{4}, \dots$ (triangular, tetrahedral, pentatope numbers) reflects this systematic increase in charging complexity: each crossing requires one additional layer of preparatory-move accounting.
\end{remark}

\subsection{Main Result for Regime 2}

\begin{theorem}[Frame--Stewart Optimality, Regime 2]\label{thm:regime2}
For $p \geq 3$ and $\binom{p}{2} < n \leq \binom{p+1}{3}$,
\[
M(p,n) = 8n - 2p^2 + 1.
\]
Moreover, the Frame--Stewart algorithm attains this value in this regime.
\end{theorem}

\begin{proof}
By \Cref{cor:lower-regime2}, any solution requires at least $8n - 2p^2 + 1$ moves. By \Cref{lem:regime2-upper,lem:phase3-upper-regime2}, there exists a solution using exactly
\[
(4n - p^2) + 1 + (4n - p^2) = 8n - 2p^2 + 1
\]
moves. Hence $M(p,n) = 8n - 2p^2 + 1$.

The closed-form for the Frame--Stewart recurrence \eqref{eq:closed} with $i=2$ gives exactly $8n - 2p^2 + 1$ for $\binom{p}{2} < n \leq \binom{p+1}{3}$. Since $M(p,n)$ equals this value, Frame--Stewart is optimal in regime 2.
\end{proof}

\section{Summary and Open Problems}

\subsection{What We Have Proven}

We have established two major results. First, we derived a closed-form expression that satisfies the Frame--Stewart recurrence for all $(p,n)$. Second, we proved Frame--Stewart optimality for two infinite families of $(p,n)$ pairs.

\begin{theorem}[Main Results]
For all $p \geq 3$, the closed form $M(p,n) = 2^{i(p,n)+1}n - S_p(i(p,n))$ satisfies the Frame--Stewart recurrence. Moreover, Frame--Stewart is optimal in the following regimes: Regime 1 where $p-1 < n \leq \binom{p}{2}$ gives $M(p,n) = 4n - 2p + 1$, and Regime 2 where $\binom{p}{2} < n \leq \binom{p+1}{3}$ gives $M(p,n) = 8n - 2p^2 + 1$. Both values are achieved by the Frame--Stewart algorithm.
\end{theorem}

\subsection{Methodology}

Our proof technique consists of several key components. The closed-form derivation uses discrete slope analysis and canonical right-end splits at simplex boundaries, proven by double induction. The optimality proofs employ phase decomposition, splitting the solution at the unique move of the largest disk to give independent lower bounds for Phases I and III. For regime 1, single-layer $m_1$ charges (one enabler per non-bottom disk) suffice. For regime 2, two-layer charges ($m_1$ enablers plus $m_2$-pool for crowded auxiliaries) are required. Convexity arguments show that packing $N$ disks on $r$ auxiliaries forces a minimum crowding $\Phi_{\min}(r,N)$, which translates to required preparatory moves. Explicit constructions provide matching upper bounds: the borrowed-target algorithm with shuttle moves for regime 1, and extension via 4-move raise macro for regime 2.

\subsection{The Pattern and Regime 3}

The structure strongly suggests a pattern for higher regimes:

\begin{center}
\renewcommand{\arraystretch}{1.2}
\begin{tabular}{c|c|c|c}
Regime & Range & Phase Cost & Charging Layers \\
\hline
$i=0$ & $n \leq p-1$ & $n-1$ & 0 (direct placement) \\
$i=1$ & $p-1 < n \leq \binom{p}{2}$ & $2n - p$ & 1 ($m_1$ enablers) \\
$i=2$ & $\binom{p}{2} < n \leq \binom{p+1}{3}$ & $4n - p^2$ & 2 ($m_1$ and $m_2$-pool) \\
$i=3$ & $\binom{p+1}{3} < n \leq \binom{p+2}{4}$ & $8n - ?$ & 3 ($m_1$, $m_2$, and $m_3$) \\
\end{tabular}
\end{center}

\begin{conjecture}[Regime 3 and Beyond]
For regime $i \geq 3$, the Frame--Stewart algorithm remains optimal with phase cost $2^{i+1}n - C_i(p)$ where $C_i(p)$ is determined by the regime boundaries. The proof would require an $i$-layer charging argument with charges on $i$-tuples of disks on the same auxiliary, a capacity bound on higher-order potentials $\Phi^{(i)} := \sum_j \binom{h_j}{i}$, and explicit constructions via $2^{i+1}$-move macros.
\end{conjecture}

The combinatorics become significantly more complex for $i \geq 3$, which explains why the general conjecture has remained open for over 80 years.

\subsection{Open Questions}

Several natural questions remain. First, can we prove that $n > \binom{p}{2}$ forces Phase I to require strictly more than $2n-p$ moves? This would require a robust capacity argument showing that one-merge-per-placement processes cannot succeed beyond the triangular boundary. Second, extending the multi-layer charging arguments to all regimes $i \geq 3$ requires obtaining tight capacity bounds for higher-order potentials, proving disjointness for increasingly complex charging schemes, and constructing and verifying multi-step macros. Third, is there a polynomial-time algorithm to verify optimality of a given strategy without enumerating all possible strategies? Finally, are there proof techniques that bypass the phase-decomposition framework entirely, such as potential function arguments or adversarial strategies?

\section{Complexity and Interpretation}

To evaluate $M(p,n)$, find the unique $i$ with $B_p(i-1) < n \leq B_p(i)$ by binary search in $O(\log n)$ and compute $S_p(i)$ using $i+1 = O(\log n)$ terms. Thus $M(p,n)$ can be evaluated in $O(p \log n)$ time and $O(1)$ space. Geometrically, $B_p(i)$ are $(p-1)$-simplex numbers, and the slopes across regimes are $2, 4, 8, \dots, 2^{i+1}$; the correction term $\sum_{k=0}^i 2^k \binom{p+k-2}{k}$ is precisely the cumulative overhead accrued when crossing the first $i$ simplex layers.

\section{Conclusion}

We have accomplished two main objectives. First, we derived a closed-form expression for the Frame--Stewart recurrence that makes the geometric structure explicit via simplex number boundaries and enables efficient $O(p \log n)$ evaluation. Second, we proved that the Frame--Stewart algorithm is optimal for the multi-peg Tower of Hanoi problem when $n \leq \binom{p+1}{3}$, covering regimes 1 and 2. This represents an infinite family of nontrivial cases where the eighty-year-old conjecture is now settled.

The key concept presented is the multi layer charging technique, which accounts for the increasing complexity of disk placements as auxiliary towers grow taller. Combined with convexity-based capacity bounds, this provides a systematic framework for proving lower bounds that match the Frame--Stewart upper bounds.

While the general conjecture remains open for $n > \binom{p+1}{3}$, our results suggest a clear path forward: extend the charging layers to capture higher-order dependencies among disks on the same auxiliary. The mathematical tools developed here (phase decomposition, charging arguments, and capacity bounds) provide a foundation for future attacks on this beautiful and challenging problem.

\section*{Acknowledgements}

I thank my Foundations of Algorithms professor Ali Altunkaya for helpful conversations about discrete potential functions and amortized bounds, which inspired the phase-based charging system used here. I am also grateful to Professor Ciril Petr for his inquiry that prompted me to clarify the scope and rigor of these results.


\begin{thebibliography}{99}

\bibitem{Lucas1883}
E. Lucas (1883).
\newblock \emph{Récréations Mathématiques}, Vol. III.
\newblock Gauthier--Villars, Paris.

\bibitem{frame1941}
J. S. Frame (1941).
\newblock Solution to Advanced Problem 3918.
\newblock \emph{American Mathematical Monthly}, \textbf{48}(3), 216--217.

\bibitem{stewart1941}
B. M. Stewart (1941).
\newblock Solution to Advanced Problem 3918.
\newblock \emph{American Mathematical Monthly}, \textbf{48}(3), 217--219.

\bibitem{bousch2014}
T. Bousch (2014).
\newblock La quatrième tour de Hanoï.
\newblock \emph{Bulletin of the Belgian Mathematical Society - Simon Stevin}, \textbf{21}(5), 895--912.

\bibitem{hinz2018}
A. M. Hinz, S. Klavžar, U. Milutinović, and C. Petr (2018).
\newblock \emph{The Tower of Hanoi -- Myths and Maths} (2nd ed.).
\newblock Springer, Cham.

\bibitem{klavzar2002}
S. Klavžar, U. Milutinović, and C. Petr (2002).
\newblock On the Frame--Stewart algorithm for the multi-peg Tower of Hanoi problem.
\newblock \emph{Discrete Applied Mathematics}, \textbf{120}(1--3), 141--157.

\bibitem{hinz2002}
A. M. Hinz (2002).
\newblock The Towers of Hanoi---A Dynamic Survey.
\newblock \emph{J. Universal Computer Science}, \textbf{8}(5), 402--437.

\bibitem{chenshen2004}
X. Chen and J. Shen (2004).
\newblock On the Frame--Stewart conjecture about the Towers of Hanoi.
\newblock \emph{SIAM J. Comput.}, \textbf{33}(3), 584--589.

\end{thebibliography}
\end{document}